\documentclass{elsart-preprint}

\usepackage{amssymb}
\usepackage{graphicx}
\usepackage[latin1]{inputenc}

\newcommand{\beq}{\begin{eqnarray}}
\newcommand{\eeq}{\end{eqnarray}}
\newcommand{\bq}{\begin{equation}}
\newcommand{\eq}{\end{equation}}

\newcommand{\eps}{{\eps}}
\newcommand{\N}{\mathbb N}

\newcommand{\R}{{\mathbb R}}
\newcommand{\Sphere}{{\mathbb S^2}}
\newcommand{\Z}{\mathbb Z}
\newcommand{\C}{\mathbb C}
\def\1{\mathbb I}
\renewcommand{\(}{\left(}
\renewcommand{\)}{\right)}
\renewcommand{\eps}{\varepsilon}

\newtheorem{theorem}{Theorem}
\newtheorem{lemma}[theorem]{Lemma}

\newtheorem{corollary}[theorem]{Corollary}
\newtheorem{proposition}[theorem]{Proposition}

\newcommand{\proof}{\par\smallskip\noindent{\sl Proof. \/}}
\newcommand{\finprf}{\unskip\null\hfill$\;\square$\vskip 0.3cm}

\newcommand{\eqref}[1]{(\ref{#1})}
\newcommand{\seq}[1]{(#1_n)_{n\in\N}}
\newcommand{\mass}{\alpha}
\newcommand{\masscor}{\beta}

\begin{document}

\begin{frontmatter}

\title{Multiplicity results for the assigned Gauss curvature problem in $\R^2$}

\author[Ceremade]{Jean Dolbeault}
\ead{dolbeaul@ceremade.dauphine.fr}
\ead[url]{www.ceremade.dauphine.fr/$\sim$dolbeaul}
\author[Ceremade]{Maria J. Esteban}
\ead{esteban@ceremade.dauphine.fr}
\ead[url]{www.ceremade.dauphine.fr/$\sim$esteban}
\author[TorVergata]{Gabriella Tarantello}
\ead{tarantel@mat.uniroma2.it}

\address[Ceremade]{Ceremade (UMR CNRS no. 7534), Universit\'e Paris-Dauphine, Place de Lattre de Tassigny, 75775 Paris Cedex~16, France.}

\address [TorVergata]{Dipartimento di Matematica. Universit\`a di Roma ``Tor Vergata", Via della Ricerca Scientifica, 00133 Roma, Italy.}

\begin{abstract}
\small To study the problem of the assigned Gauss curvature with conical singularities on Riemanian manifolds, we consider the Liouville equation with a single Dirac measure on the two-dimensional sphere. By a stereographic projection, we reduce the problem to a Liouville equation on the euclidean plane. We prove new multiplicity results for bounded radial solutions, which improve on earlier results of \hbox{C.-S.} Lin and his collaborators. Based on numerical computations, we also present various conjectures on the number of unbounded solutions. Using symmetries, some multiplicity results for non radial solutions are also stated.
\end{abstract}

\end{frontmatter}

\section{Introduction}\label{intro}

In recent years much attention has been devoted to the study of mean field equations of Liouville type on Riemann surfaces and in the presence of singular sources. Such an interest has originated from various areas of mathematics and physics, starting with the \emph{assigned Gauss curvature} problem which can be reduced to analyze the solution set of
\bq\label{GC} \Delta u+K(x)\,e^{2u} = 0\;\;\mbox{in}\;\;\R^2\eq
for a given function $K$ defined in $\R^2$. If $u$ is a solution to \eqref{GC}, then the metric $g=e^{2u}\,|dx|^2$ is conformal to the flat metric $|dx|^2$ and such that $K$ is the Gaussian curvature of the new metric $g$. Equation \eqref{GC} also appears in the analysis of gravitating systems, in the statistical mechanics description of the vorticity in fluid mechanics (see \cite{CLMP1,CLMP2,CK,Ki1,Ki2}) and has been studied more recently in the context of self-dual gauge field vortices (see \cite{JT,Y,T2}).

The solution set of equation  \eqref{GC} depends very much on the properties of $K$.  When $K$ is negative, uniqueness results are always available, while for total positive curvature $K$, either uniqueness or multiplicity of solutions holds, depending on $K$. See \cite{Li1} for various examples. In this paper we will focus on particular cases of positive functions $K$.

We are also interested in the problem of the \emph{assigned Gauss curvature with conical singularities\/} (see \cite{Tr1,Tr2,LT,E,CL3}). For a given Riemann surface $(M,g)$, we aim at determining the range of the parameters $\lambda$, $\rho\in\R$ such that
\bq\label{1.1}
\Delta_g u +\lambda\left(\frac{e^{2u}}{\int_Me^{2u}\,d\sigma_g}-\frac1{|M|}\right)-2\pi\rho\left(\delta_P-\frac1{|M|}\right)=f\,,
\eq
is solvable on $M$, where $\Delta_g$ is the Laplace-Beltrami operator, $d\sigma_g$ is the volume element corresponding to the metric $g$, $f\in C(M)$ with $\int_Mf\,d\sigma_g=0$ and $\delta_P$ is the Dirac measure with singularity at $P\in M$. In case $M$ has a non-empty boundary, both Dirichlet or Neumann boundary conditions on $\partial M$ are of interest for the applications. One could also consider sums of Dirac measures located at several source points.

In the applications, the simplest situations correspond to the $2$-sphere $M=\Sphere$ and the flat $2$-torus $M=\C/(\xi_1\Z+\xi_2\Z)$, with periodic cell domain generated by $\xi_1$ and $\xi_2$. Recall that the $2$-sphere, with the standard metric induced by Lebesgue's measure in $\R^3$, has already played a special role in the assigned curvature problem (see \cite{K,KW1,KW2}), while the torus is important since many vortex-like configurations naturally develop into periodic lattices.

For the sphere, $\rho=0$ corresponds to a particular case of the so-called Onsager vortex problem (see \cite{Li1,Li2}). In fact, for closed surfaces and when there is no singularity (\emph{i.e.\/} $\partial M=\emptyset$ and $\rho=0$), the solvability of \eqref{1.1} is quite well understood in terms of the topological properties of $M$. Starting with the work of Y. Y. Li in \cite{L}, subsequently completed by C.-C. Chen and C.-S. Lin in \cite{ChL1,ChL2}, we know that, when $\rho=0$, the solutions of \eqref{1.1} are uniformly bounded for any fixed $\lambda\in \R\setminus 4\pi\N$, and the Leray-Schauder degree $d_\lambda$ of the corresponding Fredholm operator can be explicitly computed: for $\lambda \in (4\pi(m-1), 4\pi m)$,~$m\in\N^*$,
\[\label{deg}
d_\lambda = 1\;\mbox{if}\ m=1\,,\quad d_\lambda=\frac{(-\chi(M)+1)\;\cdots\;(-\chi(M)+m-1)}{(m-1)!}\;\mbox{if}\;m\geq 2\,,
\]
where $\chi(M)$ is the Euler characteristics of $M$ (see \cite{ChL2,M}). Actually, for the flat $2$-torus we have that $d_\lambda=1$ also when $\lambda\in 4\pi \N$, and so \eqref{1.1} has a solution for every $\lambda\in \R$, if $\rho=0$ (see \cite{ChL2}). For the standard $2$-sphere, we have that $d_\lambda=0$ for all $\lambda>8\pi$. But, by a more precise topological argument (see \cite{Dj}), it can be shown that in this case equation~\eqref{1.1} admits a solution for any $\lambda\in \R\setminus 4\pi\N$, if $\rho=0$, and also some multiplicity results can be proved.

The situation is much more complex in presence of a Dirac measure. In fact, an expression for the Leray-Schauder degree, given by $d_\lambda=-\chi(M)+2$, is available only when $\rho\geq 1$ and $\lambda\in(4\pi,8\pi)$ (see \cite{CLW}). Hence, for $\Sphere$, such a degree formula yields no information about the solvability of \eqref{1.1}, and indeed this issue turns out to be very delicate (see \cite{Tr1,Tr2}). A similarly delicate situation occurs for the flat $2$-torus when $\lambda=4\pi$ and $\rho=2$. Consider for instance the equation
$$\Delta_g u + 4\pi\left(\frac{e^{2u}}{\int_Me^{2u}\,d\sigma_g}-\delta_P\right)=0\quad\mbox{in}\;M=\C/(\xi_1\Z+\xi_2\Z)\,,$$
with $P=0$. Then, C.-S. Lin and C. L. Wang have shown in \cite{LiW} that there is no solution in the case of a rectangular lattice (\emph{i.e.\/} $\xi_1=a$, $\xi_2=i\,b$, $a$, $b>0$), while there is a solution for a rhombus lattice (\emph{i.e.\/} $\xi_1=a$, $\xi_2=a\,e^{i\pi/3}$, $a>0$).

As usual, the cause for such surprising existence or nonexistence situations, there is a lack of compactness for the solution set. It can explained by looking at the singular Liouville equation
\bq\label{1.3}
\Delta u +e^{2u} = 2\pi\rho\,\delta_{z_0}\quad\mbox{in}\;\R^2\,.
\eq
The solutions of \eqref{1.3} are described using complex notations (with $\R^2\simeq \C$), by means of the Liouville formula:
$$u(z)=\frac 12\,\log\left(\frac{4\,|f'(z)|^2}{(1+|f(z)|^2)^2}\right)\,,$$
where $f$ is a meromorphic function. To account for the singularity at $z_0$,
we must require that it corresponds to a pole of $f$, or a zero of $f'$, or a branch point of $f$, with the appropriate order. In particular, for all solutions $u$ of \eqref{1.3} that satisfy $e^{2u}\in L^1(\R^2)$, there holds:
$$\int_{\R^2} e^{2u}\,dx=4\pi\,(1+\rho)\,,$$
(see \cite{CL1,CL2,PT}). This explains the role of $4\pi$ for the regular problem, $\rho=0$.

\medskip One of the motivations of this paper is to connect the solvability of \eqref{1.1} to some weighted Liouville-type equations in $\R^2$ which generalize \eqref{1.3}. To see how this class of equations could arise, let us first focus on problem \eqref{1.1} over the square $2$-torus. In this case we know that a solution exists if we replace the Dirac measure by a smooth function (see \cite{ChL2}). Let $\Omega=(-1,1)^2\subset\R^2$ and suppose that the source point $P$ coincides with the origin $0\in\Omega$. For given $\lambda>0$ and $\rho>-1$, denote by $u_\varepsilon $ the solution~of
\bq\label{1.4}\left\{\begin{array}{l}\displaystyle
\Delta u_\varepsilon +\lambda\left(\frac{e^{2u_\varepsilon}}{\int_\Omega e^{2u_\varepsilon}\,dx}-\frac1{|\Omega|}\right)=2\pi\rho\left(\frac{\varepsilon^2}{\pi\,(\varepsilon^2+|x|^2)^2}-\frac{c_\varepsilon}{|\Omega|}\right)\quad\mbox{in}\;\Omega\,, \cr
u_\varepsilon \quad\mbox{doubly periodic on}\;\partial\Omega\,,
\end{array}\right.
\eq
where $c_\varepsilon:=\int_\Omega \frac{\varepsilon^2}{\pi\,(\varepsilon^2+|x|^2)^2}\,dx$ converges to $1$ as $\varepsilon \to 0$. By doubly periodic on $\partial\Omega$, we mean $u(x+2e)=u(x)$ for any $x\in\partial\Omega$, with $e=(1,0)$ or $e=(0,1)$. To prove the existence of a solution of \eqref{1.1}, it is natural to investigate under which conditions on $\lambda$ and $\rho$ we can pass to the limit in \eqref{1.4}, along with a subsequence of $\varepsilon \to 0$. In other words, whenever possible, we need to establish {\sl a priori\/} estimates for $u_\varepsilon$ in suitable norms. To this end, denote by $u_{\varepsilon,0}$ the unique solution to the problem
$$\left\{\begin{array}{l} \displaystyle
\Delta u=2\pi\rho\left(\frac{\varepsilon^2}{\pi\,(\varepsilon^2+|x|^2)^2}-\frac{c_\varepsilon}{|\Omega|}\right)\quad\mbox{in}\;\Omega\,,\cr
u \quad\mbox{doubly periodic on}\;\partial\Omega\,,\quad \int_\Omega u\;dx=0\,,
\end{array}\right.
$$
which takes the form $u_{\varepsilon,0}(x)=\frac\rho{2}\,\log(\varepsilon^2+|x|^2)+\psi_\varepsilon(x)$, for some suitable function $\psi_\varepsilon$, which is uniformly bounded in $C^{2, \alpha}$-norm, with respect to $\varepsilon>0$. Then $u_\varepsilon = v_\varepsilon + u_{\varepsilon,0}$ is a solution of \eqref{1.4} if and only if $v_\varepsilon$ satisfies:
\[\label{1.6}\left\{\begin{array}{l} \displaystyle
\Delta v_\varepsilon +\lambda\left(\frac{e^{2(u_{\varepsilon,0}+v_\varepsilon)}}{\int_\Omega e^{2(u_{\varepsilon,0}+v_\varepsilon)}\,dx}-\frac1{|\Omega|}\right)=0 \quad\mbox{in}\;\Omega\,,\cr
v_\varepsilon \quad\mbox{doubly periodic on}\;\partial\Omega\,.
\end{array}\right.
\]
The function $e^{u_{\varepsilon,0}}$ is bounded from above and from below away from zero in $C^0_{\mbox{\small loc}}(\Omega\setminus\{0\})$, uniformly in $\varepsilon\in (0,1)$. Therefore, by well known estimates based on blow-up analysis (see, \emph{e.g.,\/} \cite{L,ChL1,ChL2,BT,T1}), we know that for $\lambda\not\in 4\pi\N$, then $u_\varepsilon$ is bounded uniformly in $C^{2, \alpha}_{\mbox{\small loc}}(\Omega\setminus\{0\})$ with respect to $\varepsilon>0$.

Next, we have to investigate what may happen in a neighborhood of the origin $0\in \Omega$. To this purpose, define $w_\varepsilon := v_\varepsilon +\frac 12\,\log \lambda -\frac 12\,\log \left(\int_\Omega e^{2u_\varepsilon}\,dx\right)$ and $W_\varepsilon :=e^{2\psi_\varepsilon}$ in $B_r(0)\subset \Omega$, for $r>0$ small enough. There holds:
$$\left\{\begin{array}{l}
-\Delta w_\varepsilon =(\varepsilon^2+|x|^2)^\rho\,W_\varepsilon (x)\,e^{2w_{\varepsilon}}-\frac\lambda{|\Omega|} \quad\mbox{in}\;B_r(0)\,,\vspace*{6pt}\cr\displaystyle
\int_{B_r(0)}(\varepsilon^2+|x|^2)^\rho\,W_\varepsilon (x)\,e^{2w_{\varepsilon}}\,dx = \lambda\,\frac{\int_{B_r(0)}e^{2u_\varepsilon}\,dx}{\int_\Omega e^{2u_\varepsilon}\,dx}\le\lambda\,.
\end{array}\right.
$$

If $w_\varepsilon$ was not uniformly bounded in $B_r(0)$, then we could find a sequence $\seq{\varepsilon}$ with $\lim_{n\to\infty}\varepsilon_n=0$, and a sequence $\seq x$ of points in $B_r(0)$ with $\lim_{n\to\infty}x_n=0$ such that
\bq\label{1.7}
w_{\varepsilon_n}(x_n)=\max_{\bar B_r(0)}w_{\varepsilon_n}\longrightarrow +\infty\quad\mbox{ as}\;n\to\infty\,.
\eq
For any $n\in\N$, define $s_n:=\max\big\{\varepsilon_n, |x_n|, \exp(-\frac{w_{\varepsilon_n}(x_n)}{2(1+\rho)})\big\}$ and observe that \hbox{$\displaystyle\lim_{n\to\infty}s_n=0$}. Let $R_n:= W_{\varepsilon_n}(x_n+s_n\,x)$, $B_n:= B_{r/s_n}(0)$. For $n$ large, $U_n(x):= w_{\varepsilon_n}(x_n+s_n\,x) +2\,(1+\rho)\,\log s_n$ satisfies
$$\left\{\begin{array}{l}
-\Delta U_n =\left(\left|\frac{\varepsilon_n}{s_n}\right|^2+\left|\frac{x_n}{s_n}+x\right|^2\right)^\rho R_n\,e^{2U_n}+o(1) \quad\mbox{in}\;B_n\,,\vspace*{6pt}\cr\displaystyle
U_n(0)=w_{\varepsilon_n}(x_n)+2\,(1+\rho)\,\log s_n\,,\vspace*{6pt}\cr\displaystyle
\int_{B_n}\left(\left|\frac{\varepsilon_n}{s_n}\right|^2+\left|\frac{x_n}{s_n}+ x\right|^2\right)^\rho R_n\,e^{2U_n}\,dx\le\lambda\,.
\end{array}\right.
$$
By definition of $s_n$, we know that $\limsup_{n\to\infty}s_n\,\exp(\frac{w_{\varepsilon_n}(x_n)}{2(1+\rho)})\ge1$. We do not know whether this limit is finite or not. If
\bq\label{Formal}
\limsup_{n\to\infty}s_n\,\exp\(\frac{w_{\varepsilon_n}(x_n)}{2(1+\rho)}\)<\infty\,,
\eq
using Harnack's estimates, we can determine a subsequence along which
$$U_n\;\to\;U_\infty\;\mbox{ in }\;C^{2, \alpha}_{\mbox{\small {loc}}}(\R^2)\;,\quad \frac{\varepsilon_n}{s_n}\to \varepsilon_\infty\in[0,1]\;,\quad \frac{x_n}{s_n}\to x_\infty\in B(0,1)\subset \R^2$$
and $W_{\varepsilon_n}(x_n+s_n\,x)$ pointwise converges to a positive constant $W_\infty>0$. In addition, $U_\infty$ satisfies:
$$\left\{\begin{array}{l}
-\Delta U_\infty =(\varepsilon_\infty^2+|x_\infty+x|^2)^\rho\,W_\infty\,e^{2U_\infty} \quad\mbox{in}\;\R^2\,,\vspace*{6pt}\cr\displaystyle
U_\infty(0)=\max_{\R^2} U_\infty \geq 0\,,\vspace*{6pt}\cr\displaystyle
\int_{\R^2}(\varepsilon_\infty^2+|x_\infty+x|^2)^\rho\,W_\infty\,e^{2U_\infty}\,dx\le\lambda\,.
\end{array}\right.$$

If $\varepsilon_\infty=0$, then by well known classification results (see \cite{CL2,PT}), we know that
$$
W_\infty\int_{\R^2}|x_\infty+x|^{2\rho}\,e^{2U_\infty}\,dx =4\pi\,(1+\rho)\,.
$$
So, we could rule out the occurrence of \eqref{1.7} in this case by restricting the problem to $\lambda<4\pi\,(1+\rho)$.

Hence, assume that $\varepsilon_\infty>0$. Then by adding a constant to $U$, scaling and translating we arrive at the limiting problem:
\bq\label{1.9}\left\{\begin{array}{ll}
-\Delta U =(1+|x|^2)^\rho\,e^{2U} \quad\mbox{in}\;\R^2\,, \quad &{\rm (a)}\vspace*{6pt}\cr\displaystyle
\int_{\R^2}(1+|x|^2)^\rho\,e^{2U}\,dx\le\lambda\,.&{\rm (b)}
\end{array}\right.
\eq
Now, with $\varepsilon_\infty>0$, it is much harder to identify the range of values assumed by $\int_{\R^2}(1+|x|^2)^\rho\,e^{2U}dx$ for all solutions of (\ref{1.9}a). It is no longer a singleton, but a whole interval contained in $(2\pi\,(1+\rho), 4\pi\,(1+\rho))$ and containing $(4\pi\,\max(1,\rho), 4\pi\,(1+\rho))$, whose explicit range is still under investigation, even when we restrict the problem to radially symmetric solutions of (\ref{1.9}a), see \cite{CheL,Li1,Li2}. In this paper, we are going to identify necessary conditions on $\lambda$ that allow us to rule out the occurrence of~\eqref{1.7} under Condition~\eqref{Formal}. In other words, if, for some $\lambda_0$, there is no solution of \eqref{1.9} for any $\lambda\le\lambda_0$, then \eqref{1.1} with $M=\Omega$ would have a solution, provided Condition~\eqref{Formal} holds.
\medskip

On the other hand, consider on $\R^2$ the solutions of
\bq\label{1.21}\left\{\begin{array}{l}
-\Delta u =(1+|x|^2)^N\,e^{2u}\quad\mbox{in}\;\R^2\,,\vspace*{6pt}\cr\displaystyle\int_{\R^2} (1+|x|^2)^N\,e^{2u}\,dx=\lambda\,.
\end{array}\right.
\eq
This problem is equivalent to \eqref{1.1}. Let $\Sigma :\Sphere\to\R^2$ be the stereographic projection with respect to the north pole, $\mathsf N:=(0,0,1)$. To any solution $u$ of~\eqref{1.21}, we can associate
\bq\label{Eq:TransfS2R2}
v(y):=u(x)-\frac 12\,\log\(\frac\lambda{(1+|x|^2)^{N+2}}\)-\log 2\,,
\eq
for $x=\Sigma(y)$ and any $y\in\Sphere$. The function $v$ solves
\bq\label{SN}
\Delta_g v +\lambda\,\frac{e^{2v}}{\int_\Sphere e^{2v}\,d\sigma_g}=\frac{2\pi\,(N+2)}{|\Sphere|}\quad\mbox{on}\;\Sphere\setminus\{\mathsf N\}\,.
\eq
Moreover, if $u$ is a locally bounded solution of \eqref{1.21}, according to \cite{CheL}, $u(x)\sim-\frac\lambda{2\pi}\,\log|x|$ as $|x|\to\infty$, which shows that
\[
\lim_{y\in\Sphere,\,y\to\mathsf N}\frac{v(y)}{\log|y-\mathsf N|}=\frac\lambda{2\pi}-(N+2)\,,
\]
Here, by $|y-\mathsf N|$, we denote the euclidean distance from $y$ to $\mathsf N$ in $\R^3\supset\Sphere$. As a consequence,
\[
\Delta_g v= \Big(\lambda-2\pi\,(N+2)\Big)\,\delta_{\mathsf N} +O(1)
\]
as $y\to\mathsf N$. Hence $v$ is a solution of \eqref{1.1} in $S^2$ with $\rho=\frac\lambda{2\pi}-(N+2)$ and $f=0$. Here the parameter $\rho$ is not the same as in \eqref {1.9}. Notice that $v$ given by \eqref{Eq:TransfS2R2} is such that $\int_\Sphere e^{2v}\,d\sigma_g=1$. As a special case, if $\lambda=2\pi\,(N+2)$, we find that $v$ is a bounded solution of \eqref{1.1} with $\rho=0$. Viceversa, if $v$ is a solution of \eqref{1.1} with $f=0$, then
\[
u=v\circ\Sigma^{-1}-\frac 12\,\log\(\int_{\Sphere}e^{2v}\,d\sigma_g\)+\log 2+\frac 12\,\log\(\frac{\lambda}{(1+|x|^2)^{N+2}}\)
\]
is a solution of \eqref{1.21}. So, existence and multiplicity of solutions to \eqref{1.21} provide us with existence and multiplicity of solutions to \eqref{1.1} on the $2$-sphere.

From now on, we shall focus on the study of \eqref{1.21}, for all positive values of $N$. Our main result establishes the existence of many bounded, radial  solutions.
\begin{theorem}\label{main}
For all $k\geq 2$ and $N>k(k+1)-2$, there are at least $2(k-2)+2$ distinct radial solutions of \eqref{1.21} with $\lambda =2\pi\,(N+2)$, one of them being the function $u^*_N(r):= \frac12\,\log\(\frac{2(N+2)}{(1+r^2)^{N+2}}\)$ with $r=|x|$.
\end{theorem}

To conclude this introduction let us mention that in many papers devoted to other applications than assigned Gauss curvature problems, the convention is to consider equation \eqref{1.1} with $e^u$ instead of $e^{2u}$. Necessary adaptations are straightforward (replace $u$ by $\frac 12\,(u-\log 2)$) and therefore are left to the interested reader.

\section{Known results}\label{known}

In this section, we collect known existence results and some of the properties of the bounded solutions to
\bq\label{PN}
\Delta u + (1+|x|^2)^N e^{2u}=0\;,\;\;\;x\in \R^2\,,\quad\mbox{and}\;\int_{\R^2}(1+|x|^2)^N e^{2u}\,dx<\infty\,,
\eq
where $N>0$ is a real parameter which in the applications usually enters as an integer value. Several papers, mostly by C.-S. Lin \emph{et al.} \cite{CheL,Li1,Li2}, have dealt with this class of equations. Almost all results are concerned with the set of radially symmetric solutions of \eqref{PN}, which can be parametrized by a parameter $a\in \R$ as follows:
\bq\label{PNa}\left\{\begin{array}{l}
u''_a+\frac{u'_a}{r} + (1+r^2)^N e^{2u_a} =0\quad\mbox{in}\; (0, +\infty)\,,\vspace*{6pt}\cr\displaystyle
u_a(0)=a\,,\quad u'_a(0)=0\,,\quad \int_0^\infty (1+r^2)^N e^{2u_a}\,r\,dr<+\infty\,.\end{array}\right.
\eq
It is proven in \cite{CheL,Li1,Li2} that, for every $a\in\R$, there exists a unique solution~$u_a$ of \eqref{PNa}, which moreover satisfies
\bq\label{comport-asymp}
\lim_{r\to\infty}\(u_a(r)+\mass(a)\,\log r\)=\masscor(a)\,,
\eq
where $a\mapsto\mass(a)$ and $a\mapsto\masscor(a)$ are two $C^1(\R,\R)$ functions defined by:
$$\mass(a) = \int_0^\infty (1+r^2)^N e^{2u_a}\,r\,dr\quad\mbox{and}\quad \masscor(a)= \int_0^\infty (1+r^2)^N e^{2u_a}\,r\,\log r\,dr\,.$$
Moreover it was proved in \cite{CheL} (also see \cite{Li2}) that, for any $N>0$,
\bq\label{limits}
\lim_{a\to -\infty}\mass(a) = 2(N+1)\quad\mbox{and}\quad \lim_{a\to+\infty}\mass(a)=2\,\min\{1,\,N\}\,.
\eq
Pohozaev's identity applied to \eqref{PN} shows that $\mass(a)\in (2, 2(N+1))$. For integrability reasons, we also know that $\mass(a) > N+1$, and so
\[
\max\{2\,,N+1\} < \mass(a) < 2(N+1)\quad\forall\,a\in\R\,.
\]
Our problem is to find the solutions of \eqref{PNa} corresponding to a given $\mass$, that~is
\bq\label{PNaBis}\left\{\begin{array}{l}
u''+\frac{u'}{r} + (1+r^2)^N e^{2u} =0\quad\mbox{in}\; (0, +\infty)\,,\vspace*{6pt}\cr\displaystyle
u'(0)=0\,,\quad \int_0^\infty (1+r^2)^N e^{2u}\,r\,dr=\mass\,.\end{array}\right.
\eq

In \cite{Li1,Li2}, C.-S. Lin investigates the uniqueness issue for problem \eqref{PNaBis} by identifying the values of $\mass$ for which there is a unique $a$ such that $\mass(a)=\mass$. He proves in \cite{Li2} that uniqueness holds for $\mass\in (2N,2(N+1))$ if $N>1$ and for all $\mass\in(2,2(N+1))$ if $N\leq 1$. On the other hand, it is easy to verify that for all $N$, the function
\bq\label{explicit}
u^*_N(r):= \frac12\,\log\(\frac{2(N+2)}{(1+r^2)^{N+2}}\)
\eq
is a solution to \eqref{PNa} for $a=a^*_N:= \frac 12\log{({2(N+2)})}$, and satisfies $\mass(a^*_N)= N+2$. Since $N+2<2N<2(N+1)$ for all $N>2$, by continuity of $a\mapsto\mass(a)$, it appears that there exists at least two different values of $a$ such that $\mass(a)=\mass$, for any $\mass \in (\min_{a\in\R}\mass(a), 2N)$. In other words, for those values of $\mass$ there exists at least two radially symmetric solutions of \eqref{PN} satisfying:
\[\label{star}
\int_{\R^2} (1+|x|^2)^N e^{2u}\,dx = 2\pi\,\mass \,.
\]
Moreover for $N=2$, $\mass(a^*_2)=4=2N$ and $\mass'(a^*_2)<0$, so the above multiplicity results also holds true for $N=2$. Summarizing, we can state the following result.
\begin{theorem}\label{Thm}{\rm\cite{CheL,Li1,Li2}} Let $N$ be any positive real number.
\begin{enumerate}
\item[{\rm (i)}] If $N\leq 1$, then the curve $a\mapsto \mass(a)$ is monotone decreasing. Moreover, there exists a radially symmetric solution $u$ of \eqref{PNaBis} if and only if $\mass \in (2, 2(N+1))$, and such a solution is unique.
\item[{\rm (ii)}] If $N> 1$, then for all $\mass\in (2N, 2(N+1))$, there exists a unique $a\in\R$ such that $\mass(a)=\mass$. In other words, for such $\mass$, problem \eqref{PNaBis} is satisfied by a unique radial solution.
\item[{\rm (iii)}] If $N\geq 2$, then $\min_{a\in\R}\mass(a)<2N$, and for all $\mass\in (\min_{a\in\R}\mass(a),\,2N)$, there exists at least two radial solutions of \eqref{PNaBis}.
\end{enumerate}
\end{theorem}

\noindent{\bf Remark. } Concerning part (iii) of Theorem~\ref {Thm}, by a closer inspection of the results of \cite{Li1}, actually we know that for $N\ge 2$, problem \eqref{PNaBis} is satisfied by a unique radial solution also when $\mass=2N$.

An important tool in the proof of the above uniqueness results is the study of
the linearized problem
\bq\label{linearPNa}\left\{\begin{array}{l}
\varphi''_a+\frac{\varphi'_a}{r} + 2\,(1+r^2)^N e^{2u_a}\,\varphi_a =0\;,\;\;\;r\in (0, +\infty)\,,\vspace*{6pt}\cr\displaystyle
\varphi_a(0)=1\,,\quad \varphi'_a(0)=0\,,\end{array}\right.
\eq
and in particular the number of zeroes of the function $\varphi_a$ when $\mass$ is in the range $(2N,2(N+1))$. The number of critical points of $a\mapsto\mass(a)$ is also connected with the number of zeroes of $\varphi_a$ in the range $(\min_{a\in\R}\mass(a),\,2N)$. It is indeed easy to prove that as $r$ goes to $+\infty$,
\bq\label{Eq:varphiInfty}
\varphi_a(r) \sim -\mass'(a)\log r + b'(a) + o(1)\,,
\eq
and hence, that $\varphi_a$ is a bounded function if and only if $a\in\R$ is a critical point of the function $\mass$. As a special case, for all $N$, if $\min_{a\in\R}\mass(a)$ is achieved for some finite $\underline a$, then $\varphi_{\underline a}$ is bounded.


\subsection{Non radially symmetric solutions}\label{NRS-known}

By \eqref{Eq:TransfS2R2}, to any solution $u_a$ of \eqref{PNa}, we can associate a function $v_a$ on $\Sphere$, such that $\int_{\Sphere}e^{2v}\,d\sigma_g=1$, which solves \eqref{SN} for $\lambda=2\pi\,\mass(a)$. At level $\mass=N+2$, $v_a$ is a bounded solution of \eqref{1.1} (with $f=0$, $\rho=0$), which is axially symmetric with respect to the unit vector $(0,0,1)$ pointing towards the north pole $\mathsf N$ of $\Sphere$. Since $v_N^*=v_{a^*_N}$ is the unique constant solution of \eqref{1.1}, if we know the existence of more than one solution at level $\mass=N+2$, then there is an axially symmetric solution of \eqref{SN} which is not constant, and that can be thus rotated in order to be axially symmetric with respect to any vector $\mathsf e\in\Sphere\setminus\{\mathsf N, \mathsf S\}$. Let us denote by $v_{\mathsf e}$ such a solution. Applying \eqref{Eq:TransfS2R2} to $v_{\mathsf e}$, we find a solution $u_{\mathsf e}$ of \eqref{PN} which is not radially symmetric. If at level $\lambda=2\pi(N+2)$ we find $k$ solutions of \eqref{PN} different from $u^*_N$, then we get $k$ punctured spheres of non radially symmetric solutions of \eqref{PN}. Details can be found in \cite{Li2}.

\section{Multiplicity results for radially symmetric solutions}\label{new}

For given $N$, critical levels of the curve $a\mapsto\mass(a)$ determine the multiplicity of the radial solutions at a given level. The number of zeros of the solutions of the linearized problem can change only at critical points of $\mass$, see below Section~\ref{sect-lin}. In the special case $\mass=N+2$, a bifurcation argument provides us with a very precise multiplicity result, which is our main result, see Section~\ref{Sect-Main}.

\subsection{A preliminary result}\label{sect-prelim}

Let $(c_k^N)_{k=1}^{n-1}$ be the ordered sequence of all critical values, counted with multiplicity, of the curve $a\mapsto\mass(a)$, $c_0^N:=\inf_{a\in\R}\mass(a)$ and $c_n^N=2(N+1)$. Denote by $(a_k^N)_{k=1}^{n-1}$ a sequence of critical points corresponding to $(c_k^N)_{k=1}^{n-1}$ and, for any $k=1$, $2,\ldots$ $n-1$, let $\epsilon_k^N:=+2$ if $a_k^N$ is a local minimum, $\epsilon_k^N:=-2$ if $a_k^N$ is a local maximum, and $\epsilon_k^N:=0$ otherwise. Also let $\epsilon_0^N:=2$ if $\inf_{a\in\R}\mass(a)$ is not achieved, and $0$ otherwise. Let $\chi_N(\mass):=1$ if $\mass>2N$ and $0$ otherwise. The next proposition links these values with the number of solutions of \eqref{PNaBis}.
\begin{proposition}\label{Prop:Counting} Let $N$ be any positive real number. With the above notations, for any $\mass>0$, Equation \eqref{PNaBis} has exactly $\sum_{j=0}^k\epsilon_j^N-\chi(\mass)$ solutions such that $\mass(a)=\mass$ if $\mass\in(c_k^N,c_{k+1}^N)\cap(\R\setminus\{2N\})$, for any $k=0$, $1$,\ldots $n-1$.\end{proposition}
The proof is straightforward and left to the reader. We shall now focus on the study of the critical points of the curve $a\mapsto\mass(a)$. Our main tool is the linearization of \eqref{PNa}.

\subsection{Study of the linearized problem}\label{sect-lin}

In order to study the multiplicity of radial solutions for \eqref{1.21}, it is convenient to perform the Emden-Fowler transformation in the linearized equation \eqref{linearPNa}:
$$t=\log r\,, \quad w_a(t) := \varphi_a(r)\,.$$
The equation in \eqref{linearPNa} is then transformed into
\bq\label{EF} w_a''(t) + 2\,e^{2t}(1+e^{2t})^N e^{2u_a(e^t)}\,w_a(t)= 0\,, \quad t\in (-\infty, +\infty)\,.
\eq
When $a=a^*_N$, the equation for $w_N^*:=w_{a_N^*}$ reads
\bq\label{EF*} {w_N^*}''(t) + \frac{(N+2)}{2\,(\cosh t)^2}\,\,w_N^*(t)= 0\,,\quad t\in (-\infty, +\infty)\,.
\eq
With one more change of variables, $w(t)=\psi(s)$, $s=\tanh t$, we find Legendre's equation:
\[
\frac d{ds}\((1-s^2)\,\frac{d\psi}{ds}\)+\frac{N+2}2\,\psi=0
\]
which defines the Legendre polynomial of order $k\in\N^*$ if $N+2=k(k+1)$. Notice that the composition of the two above changes of variables amounts to write $s=\frac{r^2-1}{r^2+1}$. A very precise spectral analysis made in \cite{Landau-Lifschitz-67} shows that the above equation has bounded solutions if and only if there is a positive integer $k$ such that $1+2k =\sqrt{1+4(N+2)}$, that is, if and only if
\[
\mathfrak N(N):=\frac{-1+\sqrt {1+4(N+2)}}{2}
\]
is a positive integer. This is solved by $N=N_k:=k(k+1)-2$, $k\in\N^*$. Actually, we are interested only in $k\ge 2$ since we only deal with $N>0$. As a consequence, $a^*_N$ is a critical point of $\mass$ if and only if $\mathfrak N(N)$ is a positive integer: $N=4$, $10$, $18$, $\ldots$

For $N=N_k>0$, we know explicitly the solutions to \eqref{EF*}. They are the Legendre polynomials, namely, with $s=\tanh t$,
\bq\label{pol} w_{N_k}^*(t)\equiv P_k(s) \;\mbox{ for all integer }\; k\geq 2\,.
\eq
\begin{lemma}\label{LLL}{\rm \cite{Landau-Lifschitz-67}} Take $N\geq 1$. Then, there exist bounded solutions of \eqref{EF*} if and only if $\mathfrak N(N)$ is a positive integer. In such a case, $\varphi_{a^*_N}$ has exactly $\mathfrak N(N)$ zeroes in the interval $(-\infty, +\infty)$.
\end{lemma}
For all $a\in\R$, $\varphi_a$ has at least two zeroes in the interval $(0, +\infty)$. This observation is a key step in the uniqueness proofs of~\cite{Li1,Li2}. Zeroes of $\varphi_a$ will also play an important role in multiplicity results. The next observation is a standard result for linear ordinary differential equations.
\begin{lemma}\label{Lem:Interval} For any $N>0$, $a_0>0$ and $R>0$, if $\varphi_{a_0}$ has $k$ zeroes in $(0,R)$ and $\varphi_{a_0}(R)\neq 0$, then there exists an $\varepsilon>0$ such that $\varphi_a$ also has exactly $k$ zeroes in $(0,R)$ for any $a\in(a_0-\varepsilon,a_0+\varepsilon)$. \end{lemma}
\proof If the result were wrong, then we could find some $r\in(0,R)$ such that $\varphi_{a_0}(r)=\varphi_{a_0}'(r)=0$ and so we would get $\varphi_{a_0}\equiv0$. \finprf
\begin{corollary}\label{corLL}
For any $N\in[N_k, N_{k+1})$, $k\in\N$, $k\ge1$, solutions to \eqref{EF*} have exactly $k+1$ zeroes in the interval $(-\infty, +\infty)$.
\end{corollary}
\proof If we normalize all functions $w_N^*$ so that $\lim_{t\to -\infty}(w_N^*(t), {w_N^*}'(t))=(1, 0)$, then by continuity in $N$ the number of zeros changes if and only if $w_N^*$ is bounded, \emph{i.e.\/} if $\mathfrak N(N)$ is a positive integer, by Lemma \ref{LLL}. At $N_1=0$, we have: $w_0^*(t)=\tanh(t)$. A careful analysis shows that, as a function of $N$, the number of zeroes is continuous from the right and increasing. \finprf

\begin{lemma} Take $N\geq 1$ and consider $a_1$, $a_2\in\R$ such that $\mass'(a_1)=\mass'(a_1)=0$ and $\mass'(a)\neq 0$ if $a\in(a_1, a_2)$. Then, for all $a\in (a_1, a_2)$, the functions $\varphi_a$ have the same number of zeroes.
\end{lemma}
\proof The proof is based on the same arguments as the proof of Corollary \ref{corLL}. One has just to replace the continuity in $N$ by the continuity in $a$. \finprf

Lemma~\ref{Lem:Interval} means that when $a$ varies, zeroes may appear or disappear only at infinity. For a given $a$, the sign of the function
\bq\label{J}
J_N(a):= \int_0^{+\infty} (1+r^2)^N e^{2u_a}\varphi_a^3\,r\,dr
\eq
governs the dynamics of the zeroes of $\varphi_a$ at infinity as follows. Denote by $r(a):=\max\{r>0\,:\,\varphi_a(r)=0\}$ the largest zero of $\varphi_a$.
\begin{lemma}\label{Lem:Dynamics} Let $\bar a>0$ be such that, for $\zeta=\pm 1$, $\lim_{a\to\bar a,\,\zeta(a-\bar a)>0}r(a)=\infty$. Then there exists $\varepsilon>0$ such that, on $(\bar a-\varepsilon,\bar a)$ if $\zeta=-1$, on $(\bar a,\bar a+\varepsilon)$ if $\zeta=+1$,
\[
\frac{dr}{da}(a)=-\frac 4{r(a)\,|\varphi_a'(r(a))|^2}\int_0^{r(a)} (1+r^2)^N\,e^{2u_a}\,\varphi^3_a\,r\,dr
\]
and $\frac{dr}{da}\,J_N(a)<0$ if $J_N(\bar a)\neq 0$. \end{lemma}
\proof First, we choose $\varepsilon>0$ small enough so that $\int_0^{r(a)} (1+r^2)^N e^{2u_a}\varphi_a^3\,r\,dr$ and $J_N(a)$ have the same sign, if $J_N(\bar a)\neq 0$. Next, we take $b>0$, small. Multiplying the equation satisfied by $\varphi_a$ and $\varphi_{a+b}$ by $r\,\varphi_{a+b}$ and $r\,\varphi_a$ respectively, and integrating by parts in the interval $(0, r(a))$, we get
\[
2\int_0^{r(a)} (1+r^2)^N\,\big(e^{2u_a}-e^{2u_{a+b}}\big)\,\varphi_a\,\varphi_{a+b}\,r\,dr= -r(a)\,\varphi'_a(r(a))\,\varphi_{a+b}(r(a))\,.
\]
By definition of $\varphi_a$ and using the uniform continuity properties of the functions $u_a$ and $\varphi_a$ on $(0, r(a))$, we obtain: $\|u_{a+b}-u_a-b\,\varphi_a\|_{L^\infty(0, r(a))}=o(b)$ as $b\to0$, $\big(1-e^{2b\varphi_a}\big)\sim-2\,b\,\varphi_a$. Since $\varphi_{a+b}(r(a))=\varphi_a'(r(a))(r(a+b))-r(a))+o(b)$, we get, as $b\to0$,
\[
-4\int_0^{r(a)} (1+r^2)^N\,e^{2u_a}\,\varphi^3_a\,r\,dr=r(a)\,|\varphi_a'(r(a))|^2\,\frac{r(a+b))-r(a)}b + o(b)\,.
\]
\finprf

Notice that, with the notations of Lemma~\ref{Lem:Dynamics}, $\varphi_{\bar a}$ is a bounded function. As a consequence, we have the following result.
\begin{corollary} Let $\tilde a$ be a critical point of $\mass$. There exists $\varepsilon>0$, small enough, such that the following properties hold.
\begin{enumerate}
\item[{\rm (i)}] If $J_N(\tilde a)>0$ and if, for any $a\in(\tilde a-\varepsilon,\,\tilde a)$, all functions $\varphi_a$ are unbounded and have $k$ zeroes in $(0, +\infty)$, then $\varphi_{\tilde a}$ is bounded and has $k$ zeroes, and for any $a\in (\tilde a,\,\tilde a +\varepsilon)$, $\varphi_a$ is unbounded and has either $k$ or $k+1$ zeroes in $(0, +\infty)$.
\item[{\rm (ii)}] If $J_N(\tilde a)<0$ and if, for any $a\in(\tilde a-\varepsilon,\,\tilde a)$, all functions $\varphi_a$ are unbounded and have $k$ zeroes in $(0, +\infty)$, then $\varphi_{\tilde a}$ is bounded and has either $k$ or $k-1$ zeroes, and for any $a\in (\tilde a,\,\tilde a +\varepsilon)$, $\varphi_a$ has the same number of zeroes as~$\varphi_{\tilde a}$.
\end {enumerate}
\end{corollary}
As already seen in Proposition~\ref{Prop:Counting}, if $\tilde a$ is a local extremum of $\mass$, the number of zeroes changes when $a$ goes through $\tilde a$, since the sign of $\varphi_a$ at infinity also changes, by \eqref{Eq:varphiInfty}. Otherwise, if $\tilde a$ is an inflection point, the number of zeroes is constant when $a$ passes through~$\tilde a$. This explains the ambiguity in the previous result.

Actually for the particular case $a=a^*_N,\, N=N_k$ we can exactly compute the value of $J_N$ as follows:
\begin{proposition}\label{signJN}
Let us define $j(k):= J_N(a^*_{N_k})$ for any integer $k\geq 2$. Then, $j(k)=0$ if $k$ is odd, and $j(k)>0$ if $k$ even.
\end{proposition}
\proof By using \eqref{pol}, we can easily compute
\bq\label{JJk}
j(k)= \frac12k(k+1)\int_{-1}^{1}P_k(s)^3\,ds\,.
\eq
Now, when $k$ is odd, $P_k$ is also odd and so, $j(k)=0$. On the contrary, Gaunt's formula, see \cite[Identity (14), page 195]{gaunt1929tha} shows that $j(k)>0$ if $k$ is even.\finprf

\noindent{\bf Remark. } If $N=N_k,\, k$ even, $a^*_N$ is a local minimum for the fonction $\mass$. Indeed, this follows from the fact that we know that $J_N(a^*_N)>0$ and the number of zeroes of $P_k$. From this we infer that the function $w^*_N$ is positive at infinity for $N$ close to $N_k$ with $N<N_k$, and negative at infinity for $N$ close to $N_k$ with $N>N_k$. Hence, $\mass$ is decreasing to the left of $a^*_N$ and increasing afterwards.

\subsection{A multiplicity result at level $\mass=N+2$}\label{Sect-Main}

This section is devoted to the proof of Theorem~\ref{main}, that is a multiplicity result for the solutions of problem \eqref{PNaBis} at level $\mass=N+2$, which also helps to illustrate Theorem 1.10 in \cite{CheL}. As seen in the introduction this amounts to study the number of bounded solutions to \eqref{1.1} with $\rho=0$ (without singularities). We show that when $\lambda = 2\pi\,\mass$ and $\mass=N+2$ becomes large, there are more and more bounded solutions to \eqref{1.1}.

Let $u$ be a radial solution of \eqref{1.21} with $\mass = N+2$. We may reformulate this problem in terms of $f:=u-u^*_N\in{\mathcal D}^{1,2}(\R^2)$ as a solution to:
\bq\label{othereq}
\Delta f +\frac\mu{(1+|x|^2)^2}\,\big(e^{2f}-1\big)=0\;\;\mbox{in}\;\;\R^2\,,\quad\int_{\R^2}\frac{e^{2f}}{(1+|x|^2)^2}\,dx=\pi\,,
\eq
with $\mu=2(N+2)$. Solutions of \eqref{othereq} are bounded by \eqref{comport-asymp}.
Moreover, \eqref{othereq} is trivially invariant under the Kelvin transformation:
\begin{lemma}\label{Lem:Kelvin} If $f$ is a solution of \eqref{othereq}, then the function $x\mapsto f\(\frac x{|x|^2}\)$ is also a solution of \eqref{othereq}. \end{lemma}
This lemma allows us to characterize many branches of solutions of \eqref{othereq}.
\begin{theorem}\label{mainBis} The function $f\equiv 0$ is a trivial solution of \eqref{othereq} for any $\mu>0$. For any $k\geq 2$, there are two continuous half-branches, $\mathcal C_k^+$ and $\mathcal C_k^-$, of solutions $(\mu,f)$ of \eqref{othereq} bifurcating from the branch of trivial solutions at, and only at, $(\mu_k=2k(k+1),0)$. Solutions in $\mathcal C_k^\pm$ are such that $\pm f(0)>0$.

Away from the trivial solutions, all branches are disjoint, unbounded and characterized by the number of zeroes. In $\mathcal C_k^\pm$, the solutions of \eqref{othereq} have exactly $k$ zeroes. If $k$ is odd, the branch $\mathcal C^\mp_k$ is the image of $\mathcal C^\pm_k$ by the Kelvin transform. If $k$ is even, the half-branches $\mathcal C^\pm_k$ are invariant under the Kelvin transform. Finally, $\mathcal C^\pm_k$ for $k\ge 3$ and $\mathcal C_2^-$ are locally bounded in $\mu$. \end{theorem}
We divide the proof in three steps. 

\subsubsection*{Step 1: Existence of unbounded branches of solutions.} We use a bifurcation method to study the set of radial solutions of \eqref{othereq} with the bifurcation parameter $\mu$. Since branches may be multi-valued in terms of~$N$, we will reparametrize them with a parameter $s$. Classical results apply for instance in $\R^+\times ({\mathcal D}^{1,2}\cap C^2)$, see \cite{MR0288640,MR0320850}.

By Lemma~\ref{LLL}, there is local bifurcation from the trivial line $\{(\mu, 0)\}$ at the points $(\mu_k:=2(N_k+2), 0)$, and there is no other bifurcation point in this trivial branch. By the properties of the Legendre polynomials $P_k$, if we denote by ${\mathcal C}^\pm_k$ the two continuous half-branches of non-trivial solutions that meet at $(\mu_k, 0)$, one easily proves that for any $(\mu, f)\in {\mathcal C}^\pm_k$ in a neighborhood of $(\mu_k, 0)$, $f(r)$ has exactly $k$ zeroes in the interval $(0, +\infty)$. 

Actually, on $\mathcal C_k^\pm$, the number of zeroes of the solutions is constant, namely equal to $k$. For instance, let us prove it for $\mathcal C_k^+$ by smoothly parametrizing the branch as follows:
\bq\label{ParamBranches}\mu(s)\;\;\mbox{with}\;\;\mu(0)=2(N_k+2)\,,\;\;\mbox{and}\;\; f=f_s\;\;\mbox{with}\;\; f_0\equiv 0\,,\;\; s\in\R\,.
\eq
Let $\Lambda_k:=\{s\in(0,\infty)\,:\,f_s\,\mbox{admits exactly $k$ zeroes}\}$. Clearly $\Lambda_k$ is not empty. Since a solution $f\not\equiv 0$ of \eqref{othereq} cannot vanish at a point together with its derivative, the smoothness of the map $s\mapsto f_s$ ensures that $\Lambda_k$ is open. We check that $\Lambda_k$ is also closed. Indeed, for $(s_n)_{n\in\N}\in\Lambda_k^\N$ with $\lim_{n\to\infty}s_n=s_\infty\in(0,\infty)$, we see that the zeroes of $f_n=f_{s_n}$ cannot collide at a point $r_0\in\R^+$ since there, $\lim_{n\to\infty}f_n=f_{s_\infty}\not\equiv 0$ would vanish together with its derivatives. So, for a zero to appear or disappear, this would require the existence of a sequence $(r_n)_{n\in\N}$ with $f_n(r_n)=0$ and $\lim_{n\to\infty}r_n=\infty$. But then, $\tilde f(r)=f_{s_\infty}(1/r)$, which is still a solution of \eqref{othereq} by Lemma~\ref{Lem:Kelvin}, would vanish at $r=0$ together with its derivatives, and this is again impossible. In conclusion, $\Lambda_k=(0,+\infty)$.

Consequently, non-trivial branches with different $k$ cannot intersect or join two different points of bifurcation in the trivial branch. For any $k\ge 2$, the half-branches $\mathcal C_k^\pm$ are therefore unbounded and we can distinguish them as follows: if $(\mu,f)\in\mathcal C_k^+$, resp. $(\mu,f)\in\mathcal C_k^-$, then $f(0)>0$, resp. $f(0)<0$. 

\subsubsection*{Step 2: Symmetry under Kelvin transform} Branches of solutions of \eqref{othereq} have an interesting symmetry property. If $k$ is odd, the solutions in the branches $\mathcal C^\pm_k$ have an odd number of zeroes in $(0,+\infty)$ and so, they cannot be invariant under the Kelvin transform, because they take values of different sign at $0$ and near $+\infty$. Since $\mu_k=2(N_k+2)$ is a simple bifurcation point, the only possibility is that the branches $\mathcal C^\pm_k$ transform into each other through the Kelvin transform. Otherwise, there would be at least four half-branches bifurcating from $(\mu_k, 0)$, which is impossible.

If $k$ is even, the solutions of $\mathcal C^\pm_k$ have an even number of zeroes. So, they take values of the same sign at $0$ and near $+\infty$. If they were not invariant under the Kelvin transform, we would find two new branches, $\tilde{\mathcal C}^\pm_k$, bifurcating from $(\mu_k, 0)$, which is again impossible.

\subsubsection*{Step 3: Asymptotic behaviour of the branches}
Non trivial branches of radial solutions are contained in the region $\{\mu>8\}$, that is $N>2$ (see Theorem~\ref{Thm} and the remark immediately afterwards for $N=2$). Furthermore, with the notations of Theorem~\ref{Thm}, by \eqref{limits}, there exists a unique $a(N)\in\R$ such that $\mass(a)>2N$ for all $a<a(N)$. For $N>2$, since $N+2<2N<2(N+1)$, if $\mass(a)=N+2$, then $a>a(N)$. Hence $f(0)=u(0)-u_N^*(0)>a(N)-u_N^*(0)$ with $u=u_a$ given by \eqref{PNaBis}.

As a consequence, the branches $\mathcal C_k^-$ are locally bounded for $\mu\in[8,+\infty)$ for any $k\ge 2$. By Step 2, $\mathcal C_3^+$ is also locally bounded for $\mu\in[8,+\infty)$. Since non trivial branches do not intersect, $\mathcal C_k^\pm$, $k\ge 3$, are all locally bounded for $\mu\in[8,+\infty)$.
\finprf

As a simple consequence of Theorem~\ref{mainBis}, we have the following corollary, which is the counterpart of Theorem~\ref{main} written for \eqref{othereq}.
\begin{cor}\label{mainBisCor}
For all $k\geq 2$, for all $\mu>\mu_k=2k(k+1)$, there at least $2(k-2)+2$ distinct radial solutions of \eqref{othereq}, one of them being the zero solution. \end{cor}

The bifurcation diagram obtained for equation \eqref{othereq} (see Fig. 1, left) is easily transformed into a bifurcation diagram for the solutions of \eqref{1.21} with $\lambda = 2\pi(N+2)$ (see Fig. 1, right) through the transformation $u= f+u^*_N$. In the case of equation \eqref{1.21}, branches bifurcate from the set of trivial solutions $\, \mathcal C :=\{(N, \, \frac12\log(2(N+2)))$, in the representation $(N, a=u(0))$.

\begin{figure}[h]\begin{center}\includegraphics[width=6cm]{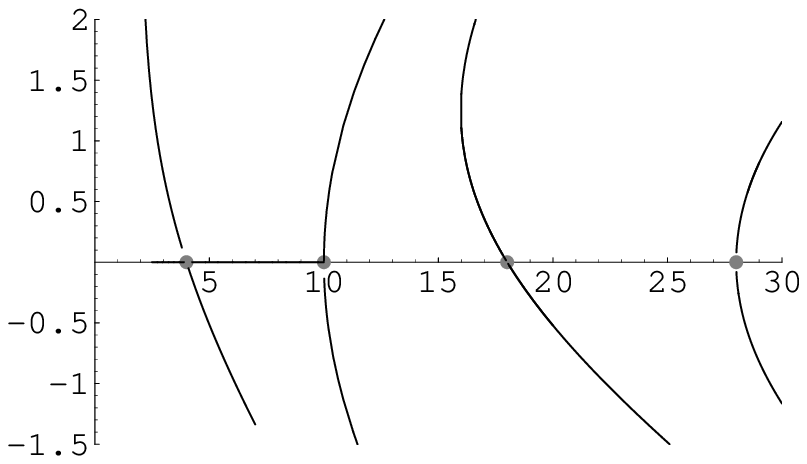}\hspace*{1.5cm}\includegraphics[width=6cm]{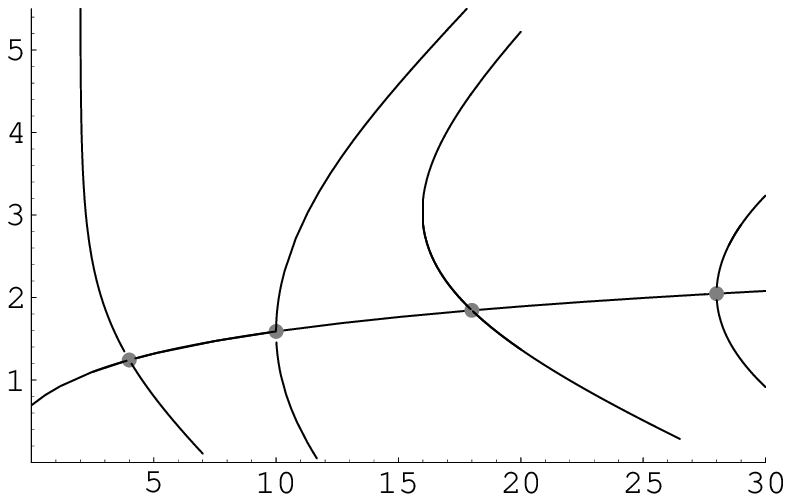}\caption{\it Bifurcation diagram in the representation $(N, f(0))$ for equation \eqref{othereq} (left) and $(N,a)$ for equation \eqref{1.21} with $\lambda = 2\pi(N+2)$ (right). Non trivial branches bifurcate from $N_k=4$, $10$, $18$, $28$,\ldots}\end{center}\end{figure}

Based on numerical evidence (see Fig. 1), it is reasonable to conjecture that, in contrast, the branch $\mathcal C_{k=2}^+$ admits a vertical asymptote in the sense that as $s\to+\infty$, then $N$ converges to $2$, which is the only admissible value by \eqref{limits}. So for $(\mu(s),f_s)\in\mathcal C_2^+$, $f_s$ should develop a concentration phenomenon at the origin, and as $s\to+\infty$, we should have: $a\to+\infty$, $N\to 2_+$ and $\frac{\mu(s)}{(1+|x|^2)^2}\,e^{2f_s}\rightharpoonup8\pi\,\delta_{z=0}$, weakly in the sense of measures.

\subsection{Non radially symmetric solutions}\label{NRS-new}

As already described in Section \ref{NRS-known}, to any solution $u\neq u^*_N$ of \eqref{PN} such that $\int_{\R^2}(1+|x|^2)^N\,e^{2u}\,dx=2\pi\,(N+2)$, we can associate a punctured sphere of non radially symmetric solutions of \eqref{PN}, $u_{\mathsf e}$ with $\mathsf e\in\Sphere\setminus\{\mathsf N, \mathsf S\}$, satisfying also $\int_{\R^2}(1+|x|^2)^N\,e^{2u_{\mathsf e}}\,dx=2\pi\,(N+2)$ for all $\mathsf e\in\Sphere\setminus\{\mathsf N, \mathsf S\}$. And so, for $N>N_k$, there are at least $2(k-2)+1$ punctured spheres of non radially symmetric solutions to \eqref{PN} at level $\lambda=2\pi(N+2)$.

\section{Further results, numerical observations and conjectures}

In the study of \eqref{PNaBis}, multiplicity results for general values of $\mass$ are difficult to deduce from Proposition~\ref{Prop:Counting} since they require a detailed analysis of the nature of each critical point: maximum, minimum, and even more in the case of an inflection point, as well as precise estimates of the corresponding critical value. Numerically, Proposition~\ref{Prop:Counting} gives straightforward results, which can be observed directly from the plots of the curve $a\mapsto\mass(a)$ for various values of $N$, see Figs.~2,~3.

Let us give some details. We consider the solution $u_a$ of \eqref{PNa} parametrized by $a=u_a(0)$. Recall that for a given $N>1$, the curve $a\mapsto\mass(a)=\int_0^\infty (1+r^2)^N e^{2u_a}\,r\,dr$ is such that $\lim_{a\to-\infty}\mass(a)=2(N+1)$, $\lim_{a\to\infty}\mass(a)=2N$, and for $N$ large enough, its range is an interval $[\mass_N,2(N+1))$ if $\mass_N$ is achieved, or $(2N,2(N+1))$ otherwise. If $N>2$, $\mass_N:=\min_{a\in\R}\mass(a)\le\mass(a_N^*)=N+2<2N$, where $a_N^*:=\frac 12\,\log(2\,(N+2))$, which provides a multiplicity result for \eqref{PNaBis} in the range $\mass\in(\mass_N,2N)\supset(N+2,2N)$.

\begin{figure}[h]\begin{center}\includegraphics[width=6cm]{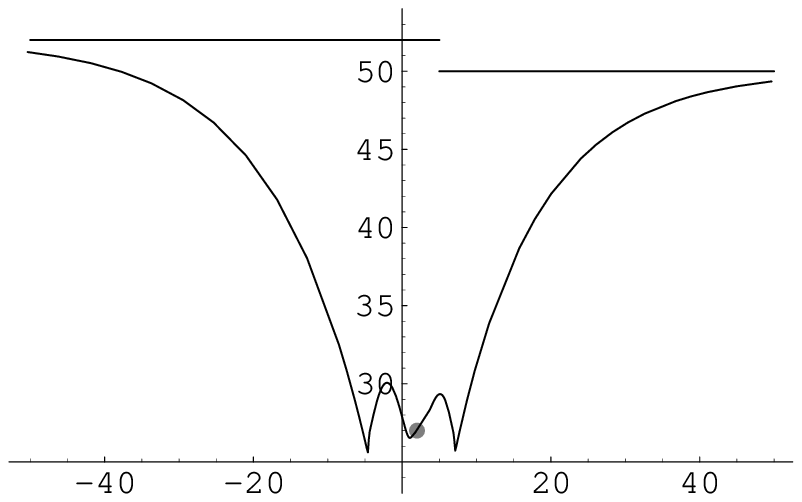}\hspace*{1.5cm}\includegraphics[width=6cm]{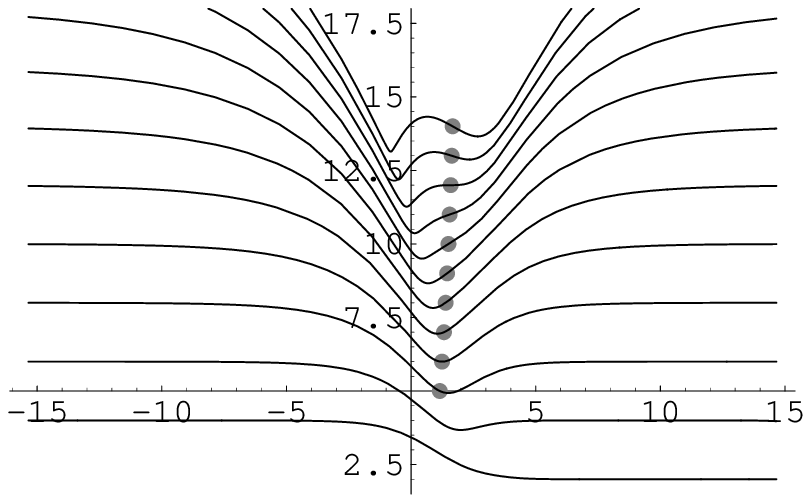}
\caption{\it Curves $a\mapsto\mass(a)$ for various values of $N$: $N=25$ (left) and $N=1$, $2$, $3$,\ldots $12$ (right). The point $(a_N^*,N+2)$ corresponding to the explicit solution \eqref{explicit} is represented by a gray dot.}
\end{center}\end{figure}

\begin{figure}[h]\begin{center}\includegraphics[width=6cm]{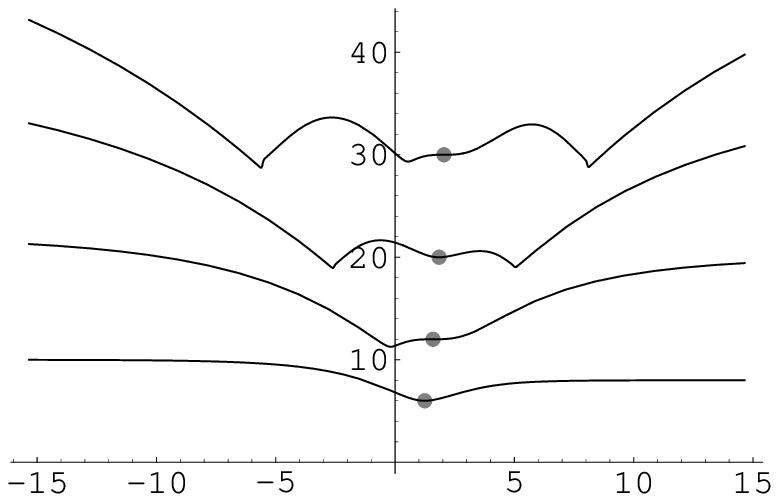}\hspace*{1.5cm}\includegraphics[width=6cm]{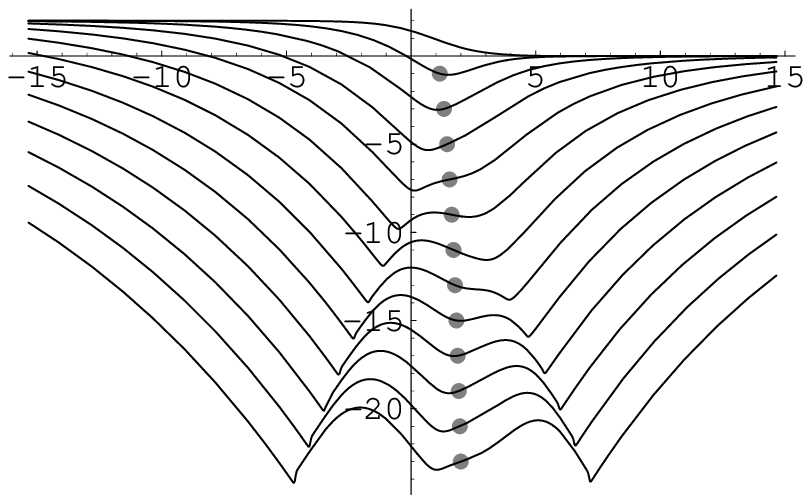}
\caption{\it Curves $a\mapsto\mass(a)$ for $N=N_k$, $k=2$, $3$, $4$, $5$ (left) and $a\mapsto\mass(a)-2N$ for $N=1$, $3$, $5$,\ldots $19$. The function $N \mapsto\mass_N-2N$ is monotone decreasing.}
\end{center}\end{figure}

As a function of $N>0$, we observe that $\mass_N=\inf_{a\in\R}\mass(a)<2N$ if and only if $N>N_0$, where $N_0$ is numerically found of the order of $1.27\pm0.02$, although its exact value is not easy to determine, see Fig. 4, left. In the range $(N_0,20)$, we observe that $\mass_N$ is achieved by the first critical point, see Fig. 4, right.

\begin{figure}[h]\begin{center}\includegraphics[width=6cm]{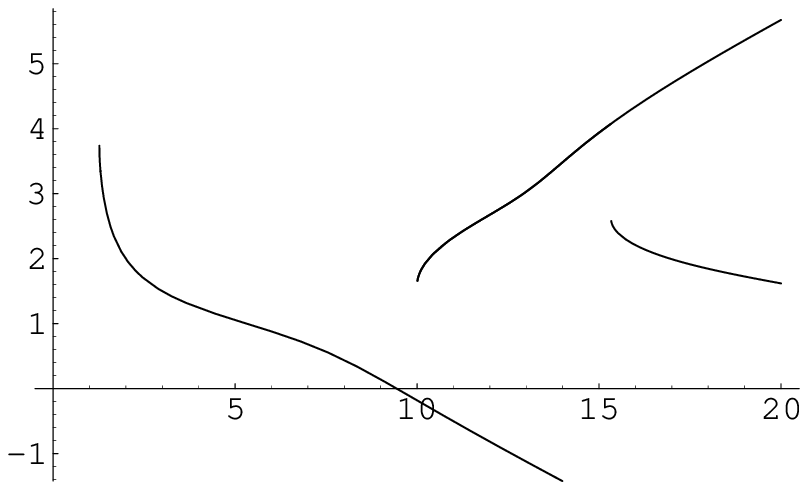}\hspace*{1.5cm}\includegraphics[width=6cm]{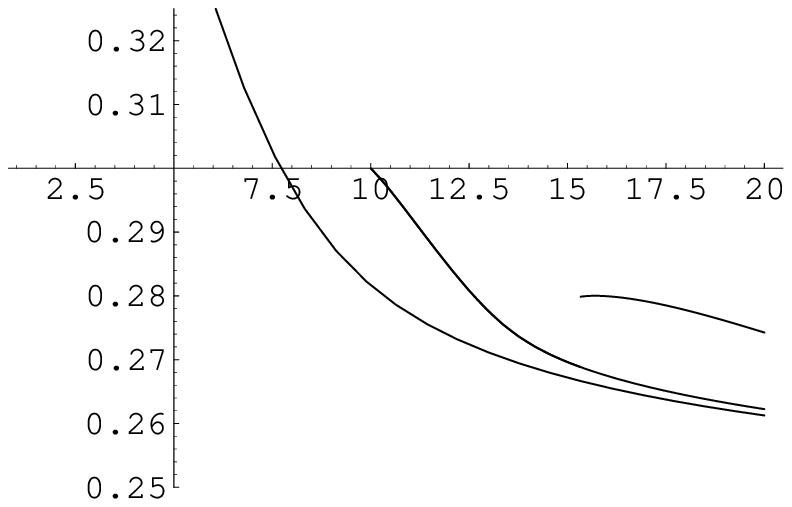}
\caption{\it Critical points (left) and critical values (divided by $4N$, right) of $\mass$, as a function of $N$.}
\end{center}\end{figure}

For any $N>N_0$, we observe that there are at least two radial solutions of~\eqref{PNaBis} for any $\mass\in (\mass_N, 2N)$. By the observations of Section~\ref{NRS-new}, it follows that there should be at least a radial solution and a whole punctured sphere of non radially symmetric solutions at level $\mass=N+2$. This supports a conjecture by C.-S. Lin in~\cite{Li2}. Actually, we can state the following result, which slightly improves on Theorem~\ref{Thm}, (iii), and rigorously defines $N_0$.
\begin{proposition}\label{Prop:Thm} There exists $N_0\in(1,2)$ such that, for all $N\in(N_0,\infty)$, $\mass_N=\min_{a\in\R}\mass(a)<2N$, and for all $\mass\in (\mass_N,\,2N)$, there exists at least two solutions of \eqref{PNaBis}.\end{proposition}
\proof The proof relies on the continuity of the curve $a\mapsto\mass(a)$ with respect to $N$ and the fact that $\mass'(a_2^*)<0$, thus proving that $\mass_N<2N$ also for $N<2$, close enough to $2$.
\finprf

It seems that as $N\to N_0$, $N>N_0$, $\mass_N$ is achieved by a unique $\tilde a_N\to\infty$, so we may conjecture that {\sl \eqref{PNaBis} has multiple solutions for $\mass\in(\mass_N,2N)$ if and only if $N\in(N_0,\infty)$, for some $N_0\in (1,2)$; for $N\in(0,N_0)$, the solution is unique, whenever it exists.} A possible way to tackle such a conjecture could be to show that for $N>1$, the function $N\mapsto\mass_N-2N$ is monotonically nonincreasing in~$N$, as it appears to be the case in our numerical study, see Fig. 4, right, and to exploit the fact that $\tilde a_N\to\infty$ as $N\to N_0$, $N>N_0$, see Fig. 4, left.

On the basis of our numerical results, we may also conjecture that for $N_0<N<10$ and $\mass\in (\mass_N, 2N)$, there exist exactly two radially symmetric solutions of \eqref{1.21}. This conjecture is supported by the bifurcation analysis of Section~\ref{Sect-Main} concerning the specific value $\mass=N+2\in(\mass_N,2N)$ for $N>2$ and $N\neq 4$. Note that for $N=4$, $\mass_N=N+2$ should hold. As $N$ increases, the curves $a\mapsto\mass(a)$ appear to have more and more critical points. Thus, for suitable values of $\mass$, the number of solutions increases as $N$ increases. We have already checked this fact in Theorem~\ref{main} for $\lambda=2\pi\,\mass$, $\mass=N+2$, but apparently it also holds for other values of $\mass$.

\medskip The last observations are concerned with the function $a\mapsto J_N(a)$ defined by \eqref{J}. It seems that such a function always takes positive values on $(-\infty, c(N))$ and negative values on $(c(N),\infty)$, for some $c(N)>0$. See Fig. 5, left. We may formulate this as the following conjecture: {\sl There exists a function $N\mapsto c(N)$ on $(0,+\infty)$ such that $J_N(a)=0$ if and only if $a=c(N)$ and $J_N(a)>0$ if and only if $a<c(N)$.} See Fig. 5, right.

\begin{figure}[h]\begin{center}\includegraphics[width=6cm]{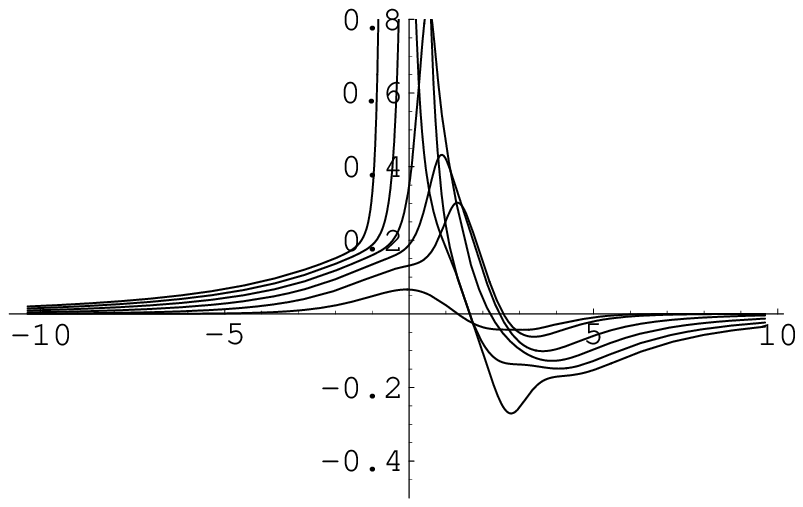}\hspace*{1.5cm}\includegraphics[width=6cm]{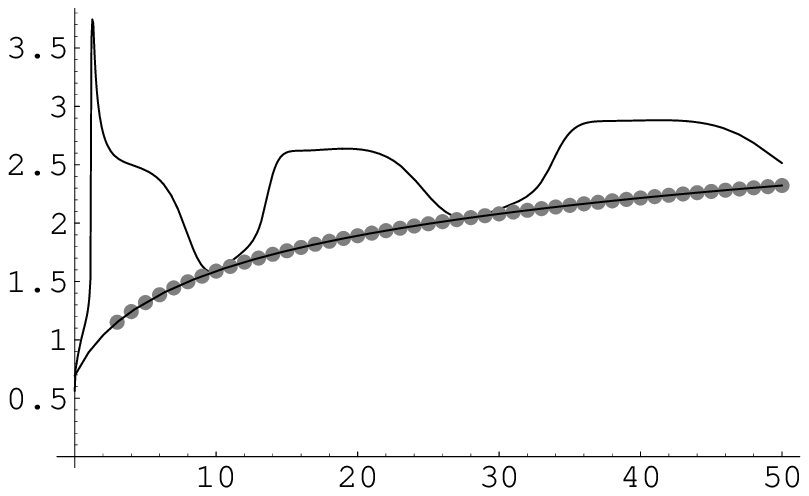}\caption{\it The function $a\mapsto J_N(a)$ for $N=1$, $3$, $5$,\ldots $11$ (left) and the curve $N\mapsto c(N)$, where, at $N$ fixed, $c(N)$ is the first positive zero of $a\mapsto J_N(a)$; the dotted line corresponds to $N\mapsto\frac 12\,\log(2(N+2))$. These two curves are tangent at $N=10=N_3$ and $N=28=N_5$ (right).}
\end{center}\end{figure}

Quite interesting is the comparison of $a_N^*$ with $c(N)$. By Proposition~ \ref{signJN}, we know that $j(k):= J_N(a^*_{N_k})=0$ if $k$ is odd and $j(k)>0$ if $k$ is even. Recall that $\mass'(a_N^*)=0$ if and only if $N=N_k$. We observe numerically, and conjecture, that: {\sl For any $N>2$, $a_N^*\le c(N)$, with equality if and only if $N=N_{2l+1}$ for some $l\ge 1$.\/} This is observed numerically with a very high accuracy for $k=3$, $5$, $7$, $9$, $11$, see Figs.~5, right, and also~6, left.

Summarizing the results of Section~\ref{sect-lin}, we have shown that, in the interval $(c(N),\infty)\ni a$, the number of nodes of $\varphi_a$ given by~\eqref{linearPNa} increases as $N$ grows each time a new critical point of $\mass$ appears. This needs to be interpreted in terms of Morse index, which is still an open question.

To investigate whether a critical point of $\mass$ is a local minimum, we may look at the functional
\[
K_N(a):=\int_0^{+\infty} (1+r^2)^N e^{2u_a}(\psi_a+2\varphi_a^2)\,r\,dr
\]
where $\psi_a$ solves the ordinary differential equation
\[\left\{\begin{array}{l}
\psi''_a+\frac{\psi'_a}{r} + 2\,(1+r^2)^N e^{2u_a}\,(\psi_a+2\,\varphi_a^2) =0\;,\;\;\;r\in (0, +\infty)\,,\vspace*{6pt}\cr\displaystyle
\psi_a(0)=0\,,\quad \psi'_a(0)=0\,.\end{array}\right.
\]
We have indeed $\mass''(a)=2\,K_N(a)$. No simple criterion for the positivity of $K_N(a)$ is known, but our numerical results at level $\mass=N+2$, see Fig. 6, right, combine very well with the results of Theorem~\ref{main} and the bifurcation diagrams shown in Fig.~1.

\begin{figure}[h]\begin{center}\includegraphics[width=6cm]{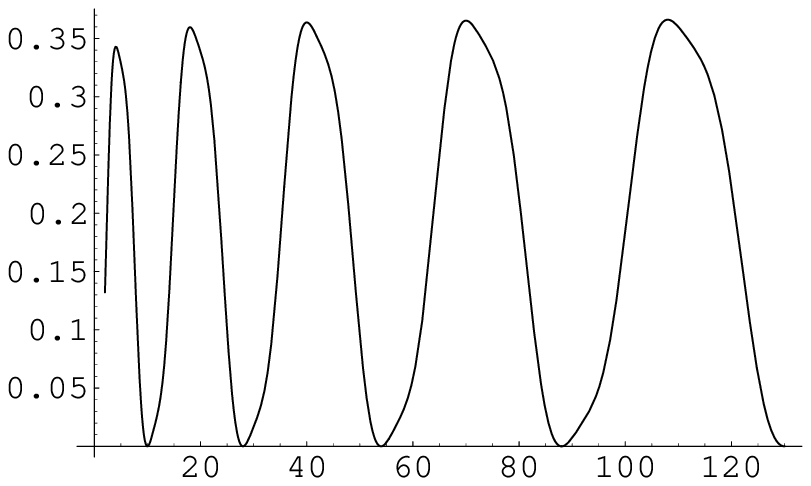}\hspace*{1.5cm}\includegraphics[width=6cm]{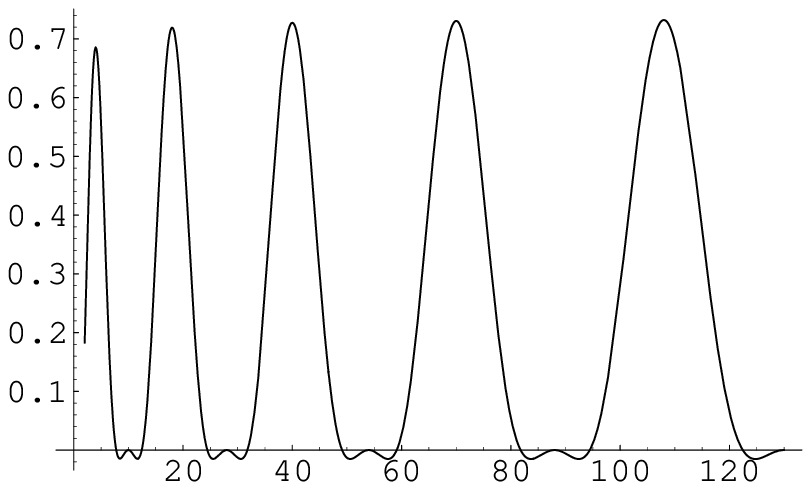}\caption{\it Left: the curve $N\mapsto J_N(a_N^*)$ is nonnegative and achieves its minimum value,~ $0$, (resp. local maxima) for $N=N_{2l+1}$, $l\ge 1$ (resp. $N=N_{2l}$). Right: the curve $N\mapsto K_N(a_N^*)$ changes sign, but is always nonnegative when $\mass'(N)=0$. When $N=N_{2l}$, $l\ge 1$, $K_N(a_N^*)$ is positive.}
\end{center}\end{figure}


\bigskip\noindent{\small{\bf Acknowlegments.} This work has been partially supported by the \emph{Fondation Sciences Mathématiques de Paris,} by the project \emph{IFO} of the French National Research Agency (ANR) and by the italian P.R.I.N. project \emph{Variational Methods and Non Linear Differential Equation,} Italy. The third author wishes also to express her gratitude to Ceremade for the warm and kind hospitality during her visits. Figures have been computed and plotted with Mathematica$^{\mbox{\tiny{\sc tm}}}$.}

\medskip\noindent{\small \copyright\,2008 by the authors. This paper may be reproduced, in its entirety, for non-commercial purposes.}


\begin{thebibliography}{10}

\bibitem{BT}
{\sc D.~Bartolucci and G.~Tarantello}, {\em Liouville type equations with
  singular data and their applications to periodic multivortices for the
  electroweak theory}, Comm. Math. Phys., 229 (2002), pp.~3--47.

\bibitem{CLMP1}
{\sc E.~Caglioti, P.-L. Lions, C.~Marchioro, and M.~Pulvirenti}, {\em A special
  class of stationary flows for two-dimensional {E}uler equations: a
  statistical mechanics description}, Comm. Math. Phys., 143 (1992),
  pp.~501--525.

\bibitem{CLMP2}
\leavevmode\vrule height 2pt depth -1.6pt width 23pt, {\em A special class of
  stationary flows for two-dimensional {E}uler equations: a statistical
  mechanics description. {II}}, Comm. Math. Phys., 174 (1995), pp.~229--260.

\bibitem{CK}
{\sc S.~Chanillo and M.~Kiessling}, {\em Rotational symmetry of solutions of
  some nonlinear problems in statistical mechanics and in geometry}, Comm.
  Math. Phys., 160 (1994), pp.~217--238.

\bibitem{ChL1}
{\sc C.-C. Chen and C.-S. Lin}, {\em Sharp estimates for solutions of
  multi-bubbles in compact {R}iemann surfaces}, Comm. Pure Appl. Math., 55
  (2002), pp.~728--771.

\bibitem{ChL2}
\leavevmode\vrule height 2pt depth -1.6pt width 23pt, {\em Topological degree
  for a mean field equation on {R}iemann surfaces}, Comm. Pure Appl. Math., 56
  (2003), pp.~1667--1727.

\bibitem{CLW}
{\sc C.-C. Chen, C.-S. Lin, and G.~Wang}, {\em Concentration phenomena of
  two-vortex solutions in a {C}hern-{S}imons model}, Ann. Sc. Norm. Super. Pisa
  Cl. Sci. (5), 3 (2004), pp.~367--397.

\bibitem{CL1}
{\sc W.~X. Chen and C.~Li}, {\em Classification of solutions of some nonlinear
  elliptic equations}, Duke Math. J., 63 (1991), pp.~615--622.

\bibitem{CL2}
\leavevmode\vrule height 2pt depth -1.6pt width 23pt, {\em Qualitative
  properties of solutions to some nonlinear elliptic equations in~{${\mathbb
  R}\sp 2$}}, Duke Math. J., 71 (1993), pp.~427--439.

\bibitem{CL3}
\leavevmode\vrule height 2pt depth -1.6pt width 23pt, {\em What kinds of
  singular surfaces can admit constant curvature ?}, Duke Math. J., 78 (1995),
  pp.~437--451.

\bibitem{CheL}
{\sc K.-S. Cheng and C.-S. Lin}, {\em On the conformal {G}aussian curvature
  equation in {$\mathbb R\sp 2$}}, J. Differential Equations, 146 (1998),
  pp.~226--250.

\bibitem{MR0288640}
{\sc M.~G. Crandall and P.~H. Rabinowitz}, {\em Bifurcation from simple
  eigen\-values}, J. Functional Analysis, 8 (1971), pp.~321--340.

\bibitem{Dj}
{\sc Z.~Djadli}, {\em Existence result for the mean field problem on {R}iemann
  surfaces of all genuses}, Commun. Contemp. Math., 10 (2008), pp.~205--220.

\bibitem{E}
{\sc A.~Eremenko}, {\em Metrics of positive curvature with conic singularities
  on the sphere}, Proc. Amer. Math. Soc., 132 (2004), pp.~3349--3355
  (electronic).

\bibitem{gaunt1929tha}
{\sc J.~A. Gaunt}, {\em The triplets of helium. {A}ppendix}, Phil. Trans. R.
  Soc. London (A), 228 (1929), pp.~192--196.

\bibitem{JT}
{\sc A.~Jaffe and C.~Taubes}, {\em Vortices and monopoles}, vol.~2 of Progress
  in Physics, Birkh\"auser Boston, Mass., 1980.
\newblock Structure of static gauge theories.

\bibitem{K}
{\sc J.~L. Kazdan}, {\em Prescribing the curvature of a {R}iemannian manifold},
  vol.~57 of CBMS Regional Conference Series in Mathematics, Published for the
  Conference Board of the Mathematical Sciences, Washington, DC, 1985.

\bibitem{KW2}
{\sc J.~L. Kazdan and F.~W. Warner}, {\em Curvature functions for open
  {$2$}-manifolds}, Ann. of Math. (2), 99 (1974), pp.~203--219.

\bibitem{KW1}
\leavevmode\vrule height 2pt depth -1.6pt width 23pt, {\em Existence and
  conformal deformation of metrics with prescribed {G}aussian and scalar
  curvatures}, Ann. of Math. (2), 101 (1975), pp.~317--331.

\bibitem{Ki1}
{\sc M.~K.-H. Kiessling}, {\em Statistical mechanics of classical particles
  with logarithmic interactions}, Comm. Pure Appl. Math., 46 (1993),
  pp.~27--56.

\bibitem{Ki2}
\leavevmode\vrule height 2pt depth -1.6pt width 23pt, {\em Statistical
  mechanics approach to some problems in conformal geometry}, Phys. A, 279
  (2000), pp.~353--368.
\newblock Statistical mechanics: from rigorous results to applications.

\bibitem{Landau-Lifschitz-67}
{\sc L.~Landau and E.~Lifschitz}, {\em Physique th\'eorique. Tome III:
  M\'ecanique quantique. {T}h\'eorie non relativiste ({F}rench)}, Deuxi\`eme
  \'edition. Translated from russian by {E}. {G}loukhian. {\'E}ditions {M}ir,
  {M}oscow, 1967.

\bibitem{L}
{\sc Y.~Y. Li}, {\em Harnack type inequality: the method of moving planes},
  Comm. Math. Phys., 200 (1999), pp.~421--444.

\bibitem{Li1}
{\sc C.-S. Lin}, {\em Uniqueness of conformal metrics with prescribed total
  curvature in~{${\mathbb R}\sp 2$}}, Calc. Var. Partial Differential
  Equations, 10 (2000), pp.~291--319.

\bibitem{Li2}
\leavevmode\vrule height 2pt depth -1.6pt width 23pt, {\em Uniqueness of
  solutions to the mean field equations for the spherical {O}nsager vortex},
  Arch. Ration. Mech. Anal., 153 (2000), pp.~153--176.

\bibitem{LiW}
{\sc C.-S. Lin and C.~L. Wang}, {\em {E}lliptic functions, {G}reen functions
  and mean field equations on torus}.
\newblock Preprint, 2006.

\bibitem{LT}
{\sc F.~Luo and G.~Tian}, {\em Liouville equation and spherical convex
  polytopes}, Proc. Amer. Math. Soc., 116 (1992), pp.~1119--1129.

\bibitem{M}
{\sc A.~Malchiodi}, {\em Morse theory and a scalar mean field equation on
  compact surfaces}.
\newblock Preprint, 2007.

\bibitem{PT}
{\sc J.~Prajapat and G.~Tarantello}, {\em On a class of elliptic problems in
  {${\mathbb R}\sp 2$}: symmetry and uniqueness results}, Proc. Roy. Soc.
  Edinburgh Sect. A, 131 (2001), pp.~967--985.

\bibitem{MR0320850}
{\sc P.~H. Rabinowitz}, {\em Some aspects of nonlinear eigenvalue problems},
  Rocky Mountain J. Math., 3 (1973), pp.~161--202.
\newblock Rocky Mountain Consortium Symposium on Nonlinear Eigenvalue Problems
  (Santa Fe, N.M., 1971).

\bibitem{T1}
{\sc G.~Tarantello}, {\em Analytical aspects of {L}iouville-type equations with
  singular sources}, in Stationary partial differential equations. Vol. I,
  Handb. Differ. Equ., North-Holland, Amsterdam, 2004, pp.~491--592.

\bibitem{T2}
\leavevmode\vrule height 2pt depth -1.6pt width 23pt, {\em Selfdual gauge field
  vortices: An analytical approach}, Progress in Nonlinear Differential
  Equations and their Applications, 72, Birkh\"auser Boston Inc., Boston, MA,
  2008.

\bibitem{Tr2}
{\sc M.~Troyanov}, {\em Metrics of constant curvature on a sphere with two
  conical singularities}, in Differential geometry (Pe\~n\'\i scola, 1988),
  vol.~1410 of Lecture Notes in Math., Springer, Berlin, 1989, pp.~296--306.

\bibitem{Tr1}
\leavevmode\vrule height 2pt depth -1.6pt width 23pt, {\em Prescribing
  curvature on compact surfaces with conical singularities}, Trans. Amer. Math.
  Soc., 324 (1991), pp.~793--821.

\bibitem{Y}
{\sc Y.~Yang}, {\em Solitons in field theory and nonlinear analysis}, Springer
  Monographs in Mathematics, Springer-Verlag New York, 2001.

\end{thebibliography}

\end{document}